\documentclass[11pt]{amsart}
\usepackage{geometry}                
\usepackage{graphicx}
\usepackage{amssymb}
\usepackage{epstopdf}
\usepackage{amsrefs}
\usepackage{enumitem}

\newtheorem{theorem}{Theorem}
\newtheorem{lemma}[theorem]{Lemma}
\newtheorem{corollary}[theorem]{Corollary}

\theoremstyle{definition}
\newtheorem{definition}[theorem]{Definition}

\theoremstyle{remark}

\newtheorem*{remarks*}{Remarks}

\newcommand\RR{\mathbb{R}}

\newcommand{\xh}{{h^x}}

\newcommand{\yh}{{h^y}}
\newcommand{\xnu}{{\nu_x}}
\newcommand{\ynu}{{\nu_y}}
\newcommand{\E}{\text{\rm e}}

\title{Non-collapsing in fully nonlinear curvature flows}
\author{Ben Andrews}
\address{Mathematical Sciences Institute, Australian National University, ACT 0200 Australia; Mathematical Sciences Center, Tsinghua University, Beijing 100084, China; Morningside Center for Mathematics, Chinese Academy of Sciences, Beijing 100190, China}
\thanks{This research was partly supported by ARC Discovery Projects grant DP0985802.  The second and third authors appreciate the support of a University of Wollongong Faculty of Informatics Research Development Scheme grant, and for the support of the Institute for Mathematics and is Applications at the University of Wollongong.}
\email{Ben.Andrews@anu.edu.au}
\author{Mat Langford}
\address{Mathematical Sciences Institute, Australian National University, ACT 0200 Australia}
\email{mathew.langford@anu.edu.au}
\author{James McCoy}
\address{Institute for Mathematics and its Applications, 
School of Mathematics and Applied Statistics, 
University of Wollongong,
Wollongong, NSW 2522, 
Australia}
\email{jamesm@uow.edu.au}

\begin{document}

\begin{abstract}
We consider embedded hypersurfaces evolving by fully nonlinear flows in which the normal speed of motion is a homogeneous degree one, concave or convex function of the principal curvatures, and prove a non-collapsing estimate:  Precisely, the function which gives the curvature of the largest interior sphere touching the hypersurface at each point is a subsolution of the linearized flow equation if the speed is concave.  If the speed is convex then there is an analogous statement for exterior spheres.  In particular, if the hypersurface moves with positive speed and the speed is concave in the principal curvatures, then the curvature of the largest touching interior sphere is bounded by a multiple of the speed as long as the solution exists.  The proof uses a maximum principle applied to a function of two points on the evolving hypersurface.  We illustrate the techniques required for dealing with such functions in a proof of the known containment principle for flows of hypersurfaces.
\end{abstract}

\maketitle

\section{Introduction}

Let $M^n$ be a compact manifold, and $X: M^n\times[0,T)\to\RR^{n+1}$ a family of smooth embeddings evolving by a curvature flow
\begin{equation}\label{E:theflow}
\frac{\partial X}{\partial t}=-F\nu,
\end{equation}
where $\nu$ is the unit normal, and the 
speed $F$ is a homogeneous degree one, monotone increasing function of the principal curvatures on a convex cone $\Gamma$ containing the positive ray.
We will assume below that $F$ is either concave or convex.  The purpose of this paper is to prove a non-collapsing result for such flows, analogous to the result proved for the mean curvature flow by the first author in \cite{ANonCollapse}.  We expect that this will provide a key step towards understanding the singular behaviour of such flows for non-convex solutions:  In the case of the mean curvature flow, the monotonicity formula of Huisken \cite{HuMono} provides a lot of information about the structure of singularities, and this is complemented by the asymptotic convexity results of Huisken and Sinestrari \cites{HS1,HS2}, and the differential Harnack or Li-Yau-Hamilton type inequality proved by Richard Hamilton \cite{HamMCFHarnack}.  The latter is available for a large class of flows \cite{AHarnack}, but there are no analogues of the monotonicity formula or the asymptotic convexity result.  The non-collapsing estimate does not precisely replace either of these, but seems nevertheless a useful tool which may be used in their stead.

The non-collapsing estimate proved for the mean curvature flow in \cite{ANonCollapse} amounts to the statement that every point of the evolving hypersurface is touched by interior or exterior spheres with radius equal to a constant $\delta$ divided by the mean curvature $H$.  It was shown there that interior non-collapsing is equivalent to the inequality 
$$
\|X(x,t)-X(y,t)\|^{2}\geq \frac{2\delta}{H(x,t)}\langle X(x,t)-X(y,t),\nu(x,t)\rangle
$$
for all $x,y\in M$.  Equivalently, this amounts to the inequality
\begin{equation}\label{eq:oldNC}
Z\left( x, y, t\right) := \frac{2\langle X(x,t)-X(y,t),\nu(x,t)\rangle}{\|X(x,t)-X(y,t)\|^{2}}\leq \frac{H(x,t)}{\delta}
\end{equation}
for all $(x,y)\in (M\times M)\setminus D$, where $D$ is the diagonal $D=\{(x,x):\ x\in M\}$.   Here we adopt the convention that the unit normal $\nu$ points outwards.
Note that the supremum of the left-hand side of \eqref{eq:oldNC} over $y$ gives the geodesic curvature of the largest interior sphere which touches at $x$.  Below we will formulate a non-collapsing result for more general curvature flows in terms of this quantity.

\begin{definition}
The \emph{interior sphere curvature} $\overline Z(x,t)$ at the point $(x,t)$ is defined by $\overline Z(x,t)=\sup\left\{Z(x,y,t):\ y\in M,\ y\neq x\right\}$.  The \emph{exterior sphere curvature} $\underline{Z}(x,t)$ at the point $(x,t)$ is defined by $\underline{Z}(x,t)=\inf\left\{Z(x,y,t):\ y\in M,\ y\neq x\right\}$.
\end{definition}

In the results to be described, an important role will be played by an equation we call the \emph{linearized flow}.  To motivate this consider a smooth family of solutions $X:\ M\times [0,T)\times (-a,a)\to\RR^{n+1}$, and define $f:\ M\times[0,T)\to\RR$ by $f(x,t) = \left\langle\frac{\partial}{\partial s}\left(X(x,t,s)\right)\Big|_{s=0},\nu(x,t)\right\rangle$.  Then $f$ satisfies the equation
\begin{equation}\label{eq:vareqn}
\frac{\partial f}{\partial t} = \dot F^{kl}\nabla_k\nabla_lf+\dot F^{kl}{h_k}^ph_{pl}f.
\end{equation}
Here $\dot F^{kl}$ is the derivative of $F$ with respect to the components $h_{kl}$ of the second fundamental form, defined by $\dot F^{kl}\big|_AB_{kl}=\frac{d}{ds}\left(F(A+sB)\right)\big|_{s=0}$ for any symmetric $B$.  Particular solutions of \eqref{eq:vareqn} include the speed $F$ (see \cite{Aconvex}*{Theorem 3.7}), corresponding to time translation $X(x,t,s) = X(x,t+s)$, the functions $\langle \nu(x,t),{\vec e}\rangle$ for $\vec{e}\in\RR^{n+1}$ fixed, corresponding to spatial translations $X(x,t,s)=X(x,t)+s\vec{e}$, and the function $\langle \nu(x,t),X(x,t)\rangle+2tF(x,t)$ (see \cite{smoczyk} or \cite{AMZConvexHypersurfaces}*{Theorem 14}), corresponding to the scalings $X(x,s,t) = (1+s)X(x,(1+s)^{-2}t)$.  

To formulate our main result we need to recall the notion of viscosity subsolution or supersolution for parabolic equations:  If $M$ is a manifold with (possibly time-dependent) connection $\nabla$ and $v:\ M\times[0,T)\to\mathbb{R}$ is continuous, then $v$ is a viscosity subsolution of the equation $\frac{\partial u}{\partial t}=G(x,t,u,\nabla u,	\nabla^{2}u)$ if for every $(x_0,t_0)\in M\times[0,T)$ and every $C^{2}$ function $\phi$ on $M\times[0,T)$ such that $\phi(x_{0},t_{0})=v(x_{0},t_{0})$, $\phi\geq v$ for $x$ in a neighbourhood of $x_{0}$ and for $t\leq t_{0}$ sufficiently close to $t_{0}$, it is true that $\frac{\partial\phi}{\partial t}\leq G(x,t,\phi,\nabla\phi,\nabla^{2}\phi)$ at the point $(x_{0},t_{0})$.  The function $v$ is a viscosity supersolution if the same holds with both inequalities for $\phi$ reversed.

Our main result is the following:  

\begin{theorem}\label{thm:main}
Assume that $X:\ M\times[0,T)\to \RR^{n+1}$ is an embedded solution of \eqref{E:theflow}.  If $F$ is convex then $\underline{Z}$ is a viscosity supersolution of the linearised flow \eqref{eq:vareqn}.
If $F$ is concave then $\overline Z$ is a viscosity subsolution of \eqref{eq:vareqn}.
\end{theorem}

Before we prove Theorem \ref{thm:main}, we mention an important consequence:

\begin{corollary}\label{cor:collapse-scale}
If $F$ is convex and positive and $X$ is an embedded solution of the curvature flow \eqref{E:theflow}, then $\inf_{M}\frac{\underline{Z}(x,t)}{F(x,t)}$ is non-decreasing in $t$.  If $F$ is concave and positive and $X$ is an embedded solution to the flow with speed $F$, then $\sup_{M}\frac{\overline Z(x,t)}{F(x,t)}$ is non-increasing in $t$.
\end{corollary}

\begin{proof}[Proof of Corollary \ref{cor:collapse-scale}]
Since $F$ satisfies equation \eqref{eq:vareqn} (see for example \cite{AMZConvexHypersurfaces}*{Lemma 9}), the result reduces to a simple comparison property of viscosity subsolutions and supersolutions.  We include the argument here for completeness:  Assume $F$ is convex, and for each $t$ let $\phi(t)=\inf_{x\in M}\frac{\underline{Z}(x,t)}{F(x,t)}$.  We must show that $\phi$ is non-decreasing in $t$.  We will accomplish this by proving that $\underline{Z}(x,t)-\left(\phi(t_0)-\varepsilon\E^{t-t_0}\right)F(x,t)\geq 0$ for any $t_0\in[0,T)$, $t\in[t_0,T)$ and $\varepsilon>0$.  Taking the limit $\varepsilon\to 0$ then gives $\underline{Z}(x,t)\geq\phi(t_{0})F(x,t)$ and hence $\phi(t)\geq\phi(t_{0})$ for $t\geq t_{0}$.

Fix $t_{0}\in[0,T)$ and $\varepsilon>0$.  Then $\underline{Z}(x,t_{0})-(\phi(t_0) -\varepsilon )F(x,t_{0})\geq \varepsilon F(x,t_{0})> 0$ for all $x$, so if $\underline{Z}-\left(\phi(t_{0})-\varepsilon\E^{t-t_{0}}\right)F$ does not remain positive for $t>t_{0}$ then there exists a time $t_{1}>t_{0}$ and a point $x_{1}\in M$ such that  $\underline{Z}-\left(\phi(t_{0})-\varepsilon\E^{t-t_{0}}\right)F$ is non-negative on $M\times[t_{0},t_{1}]$, but  $\underline{Z}(x_{1},t_{1})-\left(\phi(t_{0})-\varepsilon\E^{t_{1}-t_{0}}\right)F(x_{1},t_{1})=0$.  Since $\underline{Z}$ is a supersolution of equation \eqref{eq:vareqn}, we have at this point
\begin{align*}
0&\leq\frac{\partial}{\partial t}\left(\left(\phi(t_{0})-\varepsilon\E^{t-t_{0}}\right)F\right)-\dot F^{{kl}}\nabla_{k}\nabla_{l}\left(\left(\phi(t_{0})-\varepsilon\E^{t-t_{0}}\right)F\right)-\left(\phi(t_{0})-\varepsilon\E^{t-t_{0}}\right)F\dot F^{kl}h_{k}^{p}h_{pl}\\
&=-\varepsilon\E^{t_1-t_0}F+\left(\phi(t_{0})-\varepsilon\E^{t_1-t_{0}}\right)\left(\dot F^{kl}\nabla_k\nabla_l F+\dot F^{kl}h_k^ph_{pl}\right)\\
&\quad\null - \dot F^{{kl}}\nabla_{k}\nabla_{l}\left(\left(\phi(t_{0})-\varepsilon\E^{t_1-t_{0}}\right)F\right)-\left(\phi(t_{0})-\varepsilon\E^{t_1-t_{0}}\right)F\dot F^{kl}h_{k}^{p}h_{pl}\\
&=-\varepsilon\E^{t_1-t_0}F\\
&<0,
\end{align*}
a contradiction proving that $\underline{Z}-\left(\phi(t_{0})-\varepsilon\E^{t-t_{0}}\right)F$ remains positive.  The argument for $F$ concave is similar.
\end{proof}

Corollary \ref{cor:collapse-scale} is equivalent to the statement that the interior (for $F$ concave) or exterior (for $F$ convex) of the evolving hypersurfaces remains $\delta$-non-collapsed on the scale of $F$, in the sense of \cite{ANonCollapse}.  

We remark here that the interpretation of the non-collapsing estimate via subsolutions and supersolutions of the linearised flow \eqref{eq:vareqn} gives a new perspective even for the mean curvature flow.  Indeed, our proof is quite different from that in \cite{ANonCollapse}, and rather more transparent.  

\section{Interlude: The Containment Principle}

The proof of the main theorem uses computations of the second derivatives of the function $Z$ over the product $M\times M$, and involves a careful choice of coefficients particularly in the mixed second derivatives.   We note that there are many precedents for computations of this sort:  
Kruzhkov \cite{Kruzhkov} applied maximum principles to the difference of values at two points for solutions of parabolic equations in one space variable;  for elliptic problems quantities such as this were used by Korevaar \cite{Korevaar}, Kennington \cite{Kennington} and Kawohl \cite{Kawohl} to derive a variety of convexity properties of solutions.  For parabolic equations estimates on the modulus of continuity have been developed in \cites{AC1,AC2} and were applied in \cites{AC3,Ni} to eigenfunctions and heat kernels.  In geometric flow problems related ideas appear in work on the curve-shortening problem by Huisken \cite{HuDistComp} and Hamilton \cite{HamCSFComp} and on Ricci flow by Hamilton \cite{HamRFComp}.  More recent refinements of these techniques appear in \cites{AB1,AB2,AB3}.

Before proving the main result, we illustrate some of the techniques involved in a simpler problem:  The containment principle for solutions of fully nonlinear curvature flows of hypersurfaces.  For this problem we can consider speeds $F$ which need not be homogeneous of degree one, and need not be either convex or concave:\\

\begin{theorem}
Assume that $F$ is an odd non-decreasing symmetric function of the principal curvatures defined on $\Gamma\cup(-\Gamma)$, where $\Gamma\subset\mathbb{R}^n$ is a symmetric cone containing the positive cone, and $-\Gamma=\{-A:\ A\in\Gamma\}$.  Let $X_i: M_i \times \left[ 0, T\right) \rightarrow \mathbb{R}^{n+1}$, $i=1,2$ be two compact solutions to \eqref{E:theflow} with $X_1\left( M_1, 0\right) \cap X_2\left( M_2, 0 \right) = \emptyset$.  Then the distance from $X_1\left( M_1, t \right)$ to $X_2\left( M_2, t \right)$ is non-decreasing, and in particular  $X_1\left( M_1, t \right)\cap X_2\left( M_2, t \right)= \emptyset$ for $t\in \left[ 0, T\right)$.
\end{theorem}

\begin{proof} Define $d: M_1 \times M_2 \times \left[ 0, T\right)\to\mathbb{R}$ by
$$d\left( x, y, t\right) = \left\| X_1\left( x, t\right)- X_2\left( y,t\right) \right\| \mbox{.}$$
We show
$$\min_{M_1 \times M_2} d\left( \cdot, t\right) \geq \min_{M_1 \times M_2} d\left( \cdot, 0 \right) \mbox{,}$$
which is positive, since the initial hypersurfaces are disjoint.  As notation we will also set
$$w\left( x, y, t\right) = \frac{X_1\left( x, t\right)- X_2\left( y,t\right)}{d\left( x, y, t\right)}$$
and write $\partial_{i}^{x} = \frac{\partial X_{1}}{\partial x_{i}}$ and $\partial_{i}^{y} = \frac{\partial X_2}{\partial y_{j}}$.

The function $d$ evolves under \eqref{E:theflow} by
\begin{equation} \label{E:evlnd}
  \frac{\partial}{\partial t} d = \left< w, -F_x \nu_x + F_y \nu_y \right>.
  \end{equation}

Suppose there is a spatial minimum of $d$ at $\left( x_0, y_0, t_0 \right)$.  Then at this point,
$$\nabla^{M_1 \times M_2}d = 0 \mbox{ and } \mbox{Hess}^{M_1 \times M_2} d \geq 0 \mbox{.}$$

Choosing local orthonormal coordinates on $M_1 \times M_2$ at $\left( x_0, y_0, t_0 \right)$, that is, orthonormal coordinates $\left\{ x^i \right\}$ at $x_0$ and orthonormal coordinates $\left\{ y^i \right\}$ at $y_{0}$ we have
$$\nabla_{j}^{M_1} d = \left< \partial_{j}^{x}, w \right> \mbox{ and } \nabla_{j}^{M_2} d = - \left< \partial_{j}^{y}, w \right> \mbox{.}$$
Since we assumed that $F$ is odd, the flow is invariant under change of orientation and we can choose $\nu_x=\nu_y=w$.
In view of the definition of $w$, we have at $\left( x_0, y_0, t_0 \right)$ that
\begin{equation} \label{E:first}
\nabla_{j}^{M_1} w = \frac{1}{d} \, \partial_{j}^{x} \mbox{ and } \nabla_{j}^{M_2} w = -\frac{1}{d}\,  \partial_{j}^{y} \mbox{.}
\end{equation}

For the second spatial derivatives of $d$ we have
$$\nabla_{i}^{M_1} \nabla_{j}^{M_1} d = \left< \nabla_{i}^{M_1} \nabla_{j}^{M_1} X_1, w \right> + \left< \partial_{j}^{x}, \nabla_{i}^{M_1} w \right> \mbox{,}$$
$$\nabla_{i}^{M_2} \nabla_{j}^{M_1} d = \left< \partial_{j}^{x}, \nabla_{i}^{M_2} w \right> \mbox{ and}$$
 $$\nabla_{i}^{M_2} \nabla_{j}^{M_2} d = -\left< \nabla_{i}^{M_2} \nabla_{j}^{M_2} X_2, w \right> - \left< \partial_{j}^{y}, \nabla_{i}^{M_2} w \right> \mbox{.}$$
Using \eqref{E:first}, at $\left( x_0, y_0, t_0 \right)$ these become
$$\nabla_{i}^{M_1} \nabla_{j}^{M_1} d = \left< \nabla_{i}^{M_1} \nabla_{j}^{M_1} X_1, w \right> + \frac{1}{d}\, g_{ij}^{M_1} \mbox{,}$$
$$\nabla_{i}^{M_2} \nabla_{j}^{M_1} d = - \frac{1}{d} \left< \partial_{j}^{x}, \partial_{i}^{y} \right> \mbox{ and}$$
 $$\nabla_{i}^{M_2} \nabla_{j}^{M_2} d = -\left< \nabla_{i}^{M_2} \nabla_{j}^{M_2}, w \right> + \frac{1}{d}\, g_{ij}^{M_2}\mbox{.}$$
 
We derive the following at $\left( x_0, y_0, t_0 \right)$:  For any vector $v$ we have
\begin{align*}
0 &\; \leq v^iv^j\left( \nabla_{i}^{M_1} \nabla_{j}^{M_1} d + 2 \, \nabla_{i}^{M_2} \nabla_{j}^{M_1} d + \nabla_{i}^{M_2} \nabla_{j}^{M_2} d \right)\\
&\;=-h^x_{ij}v^iv^j\langle\nu_x,w\rangle+\frac1d g^{M_1}_{ij}v^iv^j+h^y_{ij}v^iv^j\langle\nu_y,w\rangle+\frac1d g^{M_2}_{ij}v^iv^j-\frac{2}{d}v^iv^j\langle\partial^x_i,\partial^y_j\rangle.
\end{align*}

Since $w=\nu_x=\nu_y$, the local coordinates near $x$ and $y$ may be chosen such that $\partial_{i}^{x} = \partial_{i}^{y}$ for all $i$ and $g^{M_1}_{ij}=g^{M_2}_{ij}=\delta_{ij}$.   The above becomes
$$
 h^x_{ij}v^iv^j\leq h^y_{ij}v^iv^j,
$$
or since $v$ is arbitrary, $h^x_{ij}\leq h^y_{ij}$.  Finally, since $F$ is monotone, we have $F_x\leq F_y$, and hence by \eqref{E:evlnd} we have
$$
\frac{\partial d}{\partial t}= -F_x+F_y\geq 0.
$$
\end{proof}

\begin{remarks*}
{\bf (1).} If $F$ is odd, it can be shown using a similar argument as above that for compact solutions of \eqref{E:theflow} with embedded initial hypersurface, the evolving hypersurfaces remain embedded while the curvature remains bounded.  Defining $d: M\times M \times \left[ 0, T\right)$, the curvature bound implies that there is a neigbourhood $E$ of $D= \left\{ \left( x, x\right) : x\in M \right\}$ in $M \times M$ such that 
$$d_{\RR^{n+1}}\left( x, y, \cdot \right)\geq Cd_M(x,y).$$
Consequently, the argument for the containment principle may be applied on $(M \times M)\setminus E$ to conclude that embeddedness is preserved.
 
{\bf (2).} In the containment principle the assumption that $F$ is odd can be relaxed if we make an additional topological assumption on the hypersurfaces to guarantee the correct orientation:  If we assume $F$ is defined on an arbitrary symmetric cone $\Gamma$ containing the positive cone, and $M_1=\partial\Omega_1$ and $M_2=\partial\Omega_2$ with $\Omega_1\subset\Omega_2$, and require that the unit normal to $M_i$ points out of $\Omega_i$ for $i=1,2$, then the above argument goes through with minor changes.  Without such a condition disjointness may not be preserved:  For example if $n=2$ and $F=H+|A|$, with the cone $\Gamma=\{(\kappa_1,\kappa_2):\ \max\{\kappa_1,\kappa_2\}>0\}$, then surfaces with opposite orientation having nearest points of saddle type will move closer together (and can cross).  In this example it is also true that embedded initial surfaces can evolve smoothly to become non-embedded.

\end{remarks*}
\section{Proof of the Main Theorem}

We now prove theorem \ref{thm:main}, namely, that $\underline{Z}$ ($\overline{Z}$) is a viscosity supersolution (subsolution) of the linearised flow \eqref{eq:vareqn} when if $F$ is convex (concave).

As in the previous section, the proof involves computation with the second derivatives over the product $M\times M$.   However, the computation here has an unexpected feature of the proof in the case of fully nonlinear flows:  In all the previous computations of this type mentioned above, the two points $x$ and $y$ have appeared in a symmetric way, so that the choice of coefficients in the second derivatives is determined by information at both points.  This has been a serious obstacle to applications of the methods to fully nonlinear flows, since the coefficients of the equation at different points would involve the second derivatives (or second fundamental form) at different points, and there is insufficient control on these to allow a useful comparison.   However,
in the present computation $x$ and $y$ play very different roles, and in particular the function $Z$ only depends on $x$ at the level of the highest derivatives.  Accordingly we are able to 
use a choice of coefficients in the second derivatives which depends on $x$ but not on $y$, thus removing any need to compare the second fundamental form at different points.  The key observation that makes this choice work is given in 
Lemma \ref{lem:contract-yh}.

\begin{proof}[Proof of Theorem \ref{thm:main}]
The definitions of $\overline Z(x,t)$ and $\underline{Z}(x,t)$ involve extrema of $Z$ over the noncompact set $\{y\in M:\ y\neq x\}$.  Accordingly we begin by extending $Z$ to a continuous function on a suitable compactification.  

The diagonal
$D$ is a compact submanifold of dimension and codimension $n$ in $M\times M$.   The normal subspace $N_{(x,x)}D$ of $D$ at $(x,x)$ is the subspace $\{(u,-u):\ u\in T_xM\}\subset T_{(x,x)}(M\times M)$.  The tubular neighbourhood theorem provides $r>0$ such that the exponential map is a diffeomorphism on $\{(x,x,u,-u)\in TM\times TM:\ 0<\|u\|<r\}$.
We `blow up' along $D$ to define a manifold with boundary $\hat M$ which compactifies $(M\times M)\setminus D$, as follows:
As a set, $\hat M$ is the disjoint union of $(M\times M)\setminus\{(x,x):\ x\in M\}$ with the unit sphere bundle $SM=\{(x,v)\in TM:\ \|v\|=1\}$.  The manifold-with-boundary structure is defined by the atlas generated by all charts for $(M\times M)\setminus D$, together with the charts $\hat Y$ from $SM\times(0,r)$ defined by taking a chart $Y$ for $SM$, and setting $\hat Y(z,s):=(\exp(sY(z)),\exp(-sY(z)))$.

We extend the function $Z$ to $\hat M\times [0,T)$ as follows:  For $(x,y)\in (M\times M)\setminus D$ and $t\in[0,T)$ we define
$$
Z(x,y,t) = \frac{2\langle X(x,t)-X(y,t),\nu(x,t)\rangle}{\|X(x,t)-X(y,t)\|^2}.
$$
For $(x,v)\in SM$ we define 
$$
Z(x,v,t) = h_{(x,t)}(v,v),
$$
where $h_{(x,t)}$ is the second fundamental form of $M_t$ at $x$.  Since $X$ is an embedding, $Z$ is continuous on $(M\times M)\setminus D$.  A straightforward computation shows that the above extension of $Z$ to $\hat M$ is also continuous.
It follows that $\overline Z(x,t)$ is attained on $\hat M$, in the sense that either there exists $y\in M\setminus\{x\}$ such that $\overline Z(x,t)=Z(x,y,t)$, or there exists $v\in T_xM$ with $\|v\|=1$ such that $\overline Z(x,t) = Z(x,v,t)$.  Also, since the supremum over $M\setminus\{x\}$ equals the supremum over $\hat M$, and this is no less than the supremum over the boundary $SM$, we have that $\overline Z(x,t)$ is no less than the maximum principal curvature $\kappa_{\max}(x,t)$.  Similarly, $\underline{Z}(x,t)$ is attained on $\hat M$ and is no greater than the minimum principal curvature $\kappa_{\min}(x,t)$.

To prove that $\overline Z$ is a subsolution if $F$ is concave, we consider, for an arbitrary point, $(x_0,t_0)$, an arbitrary $C^{2}$ function $\phi$ which lies above $\overline Z$ on a neighbourhood of $(x_{0},t_{0})$ in $M\times[0,t_{0}]$, with equality at $(x_{0},t_{0})$, and prove a differential inequality for $\phi$ at $(x_{0},t_{0})$.  

Observe that for all $x$ close to $x_{0}$, and all $t\leq t_{0}$ close to $t_{0}$ we have $Z(x,y,t)\leq \overline Z(x,t)\leq \phi(x,t)$ for each $y\neq x$ in $M$, and $Z(x,v,t)\leq \overline Z(x,t)\leq \phi(x,t)$ for all $v\in S_{x}M$.  Furthermore equality holds in the last inequality in both cases when $(x,t)=(x_{0},t_{0})$.  By definition of $\overline Z$ we either have $Z(x_{0},y_{0},t_{0})=\overline Z(x_{0},t_{0})$ for some $y_{0}\neq x_{0}$, or we have $Z(x_{0},\xi_{0},t_{0})=\overline Z(x_{0},t_{0})$ for some $\xi_{0}\in S_{x_{0}}M$.  

We consider the latter case first:  Define a smooth unit vector field $\xi$ near $(x_{0},t_{0})$ by choosing $\xi(x_{0},t_{0})=\xi_{0}$, extending to $(x,t_{0})$ for $x$ close to $x_{0}$ by parallel translation along geodesics, and extending in the time direction by solving $\frac{\partial\xi}{\partial t}=F{\mathcal W}(\xi)$, where $\mathcal{W}$ is the Weingarten map.  This construction implies that $\nabla\xi(x_{0},t_{0})=0$ and $\nabla^{2}\xi(x_{0},t_{0})=0$, and from the evolution equation for the second fundamental form we find that
$$
\frac{\partial}{\partial t}(h(\xi,\xi))=\dot F^{kl}\nabla_{k}\nabla_{l}(h(\xi,\xi))+\ddot F^{kl,pq}\nabla_{\xi}h_{kl}\nabla_{\xi}h_{pq}+h(\xi,\xi)\dot F^{kl}h_{k}^{p}h_{pl}
$$
at the point $(x_{0},t_{0})$.  The second term on the right is non-positive by the concavity of $F$.  At the point $(x_{0},t_{0})$ we also have $\phi=h(\xi,\xi)$, and since $\phi\geq h(\xi,\xi)$ at nearby points and earlier times we also have $\frac{\partial\phi}{\partial t}\leq\frac{\partial}{\partial t}(h(\xi,\xi))$ and $\nabla^{2}\phi\geq \nabla^{2}(h(\xi,\xi))$ at this point.  Combining these inequalities gives $\frac{\partial\phi}{\partial t}\leq \dot F^{kl}\nabla_{k}\nabla_{l}\phi+\phi\dot F^{kl}h_{k}^{p}h_{pl}$ at $(x_{0},t_{0})$ as required.

Next we consider the case where $Z(x_{0},y_{0},t_{0})=\phi(x_{0},t_{0})$ for some $y_{0}\neq x_{0}$, and $\phi(x,t)\geq Z(x,y,t)$ for all points $x$ near $x_{0}$, times $t\leq t_0$ near $t_{0}$, and arbitrary $y\neq x$ in $M$. This implies that $\frac{\partial\phi}{\partial t}(x_{0},t_{0})\leq\frac{\partial Z}{\partial t}(x_{0},y_{0},t_{0})$, that the first spatial derivatives of $\phi-Z$ in $x$ and $y$ vanish at $(x_{0},y_{0},t_{0})$ and that the second spatial derivatives of $\phi-Z$ are non-negative at $(x_{0},y_{0},t_{0})$.  We compute these derivatives, working in local normal coordinates $\{x^i\}$ near $x$ and $\{y^i\}$ near $y$. To simplify notation we define $d=|X(x,t)-X(y,t)|$ and $w=\frac{X(x,t)-X(y,t)}{d}$ and write $\partial^x_i=\frac{\partial X}{\partial x^i}$.  We first compute the first spatial derivatives with respect to $y$:
\begin{equation}\label{eq:Zy}
\frac{\partial}{\partial y^i}\left(\phi-Z\right)=\frac{2}{d^{2}}\left\langle\partial^{y}_{i},\xnu-dZw\right\rangle.
\end{equation}



The first derivatives with respect to $x$ are slightly more complicated:
\begin{equation}\label{eq:Zx}
\frac{\partial}{\partial x^{i}}\left(\phi-Z\right)=\frac{\partial\phi}{\partial x^{i}}-\frac{2}{d}\left(\xh_{i}^{p}\langle w,\partial^{x}_{p}\rangle-Z\langle w,\partial_{i}^{x}\rangle \right).
\end{equation}
The left, and therefore right, sides of equations \eqref{eq:Zy} and \eqref{eq:Zx} vanish at $(x_{0},y_{0},t_{0})$.  

Now we differentiate further to find the second derivatives: Using the fact that the first derivatives of $Z$ with respect to $y$ vanish, we find
\begin{align}
\frac{\partial^2}{\partial y^i\partial y^j}\left(\phi-Z\right)
&=\frac{2}{d^{2}}\left\{\left\langle\yh_{ij}\ynu,dZw-\xnu\right\rangle+Z\left\langle\partial^y_i,\partial^y_j\right\rangle\right\}\notag\\
&=\frac{2}{d^2}\left(Z\delta_{ij}-\yh_{ij} \right).\label{eq:Zyy}
\end{align}
Differentiating \eqref{eq:Zy} with respect to the $x$ coordinates gives the mixed partial derivatives:
\begin{equation}\label{eq:Zxy}
\frac{\partial^2}{\partial x^j\partial y^i}\left(\phi-Z \right)=-\frac{2}{d^2}\left(Z\delta_j^p-\xh_j^p\right)\langle\partial^y_i,\partial^x_p\rangle-\frac{2}{d}\frac{\partial\phi}{\partial x^j}\langle w,\partial^y_i\rangle.
\end{equation}
Differentiating \eqref{eq:Zx} with respect to the $x$ coordinates gives:
\begin{align}
\frac{\partial^2}{\partial x^i\partial x^j}\left(\phi-Z \right) &=
\frac{2}{d^2}\left(Z\delta_{ij}-\xh_{ij}\right) + Z\xh_{jp}\delta^{pq}\xh_{qi}
-\frac{2}{d}\nabla_p\xh_{ij}\delta^{pq}\langle w,\partial^x_q\rangle\label{eq:Zxx}
\\
&\quad\null -Z^2\xh_{ij}+\frac{2}{d}\frac{\partial\phi}{\partial x^j}\langle w,\partial^x_i\rangle
+\frac{2}{d}\frac{\partial\phi}{\partial x^i}\langle w,\partial^x_j\rangle+\frac{\partial^2\phi}{\partial x^i\partial x^j}.\notag\end{align}

Finally we compute the time derivative:
\begin{align}
\frac{\partial}{\partial t}\left(\phi-Z \right) &=\frac{\partial\phi}{\partial t}+\frac{2F_x}{d^2}-\frac{2F_y}{d^2}\langle \ynu,\xnu-dZw\rangle-\frac{2}{d}\langle w,\nabla F_x\rangle -Z^2F_x\notag\\
&=\frac{\partial\phi}{\partial t}+\frac{2F_x}{d^2}-\frac{2F_y}{d^2}
-\frac{2}{d}\langle w,\nabla F_x\rangle -Z^2F_x.
\label{eq:Zt}
\end{align}
Combining equations \eqref{eq:Zyy}--\eqref{eq:Zt} and the inequalities at $(x_0,y_0,t_0)$ we obtain
\begin{align}
0&\leq -\frac{\partial}{\partial t}\left(\phi-Z \right)+\dot F_x^{ij}\left(\frac{\partial^2}{\partial x^i\partial x^j}\left(\phi-Z \right)
+2\frac{\partial^2}{\partial x^i\partial y^j}\left(\phi-Z \right)+\frac{\partial^2}{\partial y^i\partial y^j}\left(\phi-Z \right)\right)\nonumber\\
&=-\frac{\partial\phi}{\partial t}+\dot F_x^{ij}\nabla_i\nabla_j\phi+\phi\dot F_x^{ij}\xh_{ip}\delta^{pq}\xh_{qj}-\frac{4F_x}{d^2}+\frac{4}{d^2}\dot F_x^{ij}\xh_{iq}\delta^{qp}\langle \partial^y_j,\partial^x_p\rangle\label{eq:Zt-LZ}\\
&\quad\null +\frac{2F_y}{d^2}-\frac{2}{d^2}\dot F_x^{ij}\yh_{ij}+\frac{4Z}{d^2}\dot F_x^{ij}\delta_{ij}-\frac{4Z}{d^2}\dot F_x^{ij}\langle \partial^x_i,\partial^y_j\rangle+\frac{4}{d}\dot F_x^{ij}\frac{\partial\phi}{\partial x^i}\langle w,\partial^x_j-\partial^y_j\rangle.\nonumber
\end{align}
Now note that, by the homogeneity of $F$, $F_x =\dot F_x^{ij}\xh_{ij}$, so that
$$
-\frac{4F_x}{d^2}+\frac{4}{d^2}\dot F_x^{ij}\xh_{iq}\delta^{qp}\langle \partial^y_j,\partial^x_p\rangle
=-\frac{4}{d^2}\dot F_x^{ij}\xh_{iq}\delta^{qp}\left(\delta_{jp}-\langle\partial^y_j,\partial^x_p\rangle\right).
$$
We can also write
$$
\frac{4Z}{d^2}\dot F_x^{ij}\delta_{ij}-\frac{4Z}{d^2}\dot F_x^{ij}\langle \partial^x_i,\partial^y_j\rangle
=\frac{4Z}{d^2}\dot F_x^{ij}\left(\delta_{ij}-\langle\partial^y_j,\partial^x_i\rangle\right).
$$
To control the first two terms on the second line of \eqref{eq:Zt-LZ} we use the following observation:

\begin{lemma}\label{lem:contract-yh}
If $F$ is concave, then for any $y\neq x$ we have
$$
\dot F_x^{ij}\yh_{ij}\geq F_y.
$$
If $F$ is convex, then the reverse inequality holds.
\end{lemma}

\begin{proof}[Proof of Lemma]
Let $A=\xh$ and $B=\yh$.  Then concavity of $F$ gives
$$
F(B)\leq F(A)+\dot F_{A}\left(B-A\right) = F(A)+\dot F_A(B)-\dot F_A(A).
$$
The homogeneity of $F$ gives by the Euler relation that $\dot F_A(A)=F(A)$, yielding
$$
F(B)\leq \dot F_A(B)
$$
as claimed.  The inequality is reversed for $F$ convex.
\end{proof}

Using these observations, together with the identity for $\frac{\partial\phi}{\partial x^i}$ coming from the vanishing of $\frac{\partial}{\partial x^i}\left(\phi-Z\right)$ in equation \eqref{eq:Zx}, we find:
\begin{align*}
0&\leq -\frac{\partial\phi}{\partial t}+\dot F_x^{ij}\nabla_i\nabla_j\phi+\phi\dot F_x^{ij}\xh_{ip}\delta^{pq}\xh_{qj}\\
&\quad\null+\frac{4}{d^2}\dot F_x^{ij}\left(Z\delta_{ip}-\xh_{ip}\right)\delta^{pq}\left(\delta_{qj}-\langle\partial^y_j,\partial^x_q\rangle +2\langle w,\partial^x_q\rangle\langle w,\partial^y_j-\partial^x_j\rangle\right).
\end{align*}

We now prove that the term in the final brackets is non-positive, that is,

\begin{lemma}
The term: $\delta_{qj}-\langle\partial^y_j,\partial^x_q\rangle +2\langle w,\partial^x_q\rangle\langle w,\partial^y_j-\partial^x_j\rangle$ is non-positive.
\end{lemma}

\begin{proof}[Proof of Lemma]
We now choose the local coordinates $\{x^i\}$ and $\{y^i\}$ more carefully. Throughout we continue to compute at the minimum $(x_0,y_0,t_0)$. Then we may choose $\partial^y_n$ and $\partial^x_n$ to be coplanar with $w$, and $\partial^y_i=\partial^x_i$ for $i=1,\dots,n-1$. This ensures that $\delta_{qj}-\langle\partial^y_j,\partial^x_q\rangle +2\langle w,\partial^x_q\rangle\langle w,\partial^y_j-\partial^x_j\rangle$ is non-zero only when $p=q=n$. 

We have two cases to consider, first supposing that $\langle w,\nu_x\rangle\geq 0$. In this case we may define $\alpha\in[0,\pi/2)$ by $\langle w,\nu_x\rangle=\sin\alpha$. Note that we have one final degree of freedom in the coordinates, namely the directions of $\partial_n^x$ and $\partial_n^y$. Direct $\partial_n^x$ such that $\langle w,\partial^x_n\rangle=-\cos\alpha$. Now define $\theta\in[0,\pi/2)$ and the orientation of $\partial^y_n$ by the conditions $\langle \partial^y_n,\partial^x_n\rangle=-\cos 2\theta$ and $\langle \partial^y_n,\nu_x\rangle=\sin 2\theta$. Then the vanishing of $\partial_{y_n}(\phi-Z)$ implies
\begin{align}
\langle \partial^y_n,\nu_x\rangle &=2\langle w,\nu_x\rangle\langle \partial^y_n,w\rangle\label{eq:ynidentity}\\
\Rightarrow \sin 2\theta\cos 2\alpha &=\sin 2\alpha\cos 2\theta\,.\nonumber
\end{align}
That is, $\sin(2\theta-2\alpha)=0$ and we find $\theta=\alpha$. The identity \eqref{eq:ynidentity} now implies that $\langle \partial^y_n,w\rangle=\cos\theta$ and we may compute,
\begin{align*}
\delta_{qj}-\langle\partial^y_j,\partial^x_q\rangle +2\langle w,\partial^x_q\rangle\langle w,\partial^y_j-\partial^x_j\rangle& 
=1+\cos(2\theta)+2\cos\theta(-\cos\theta-\cos\theta)\\
& = 2\cos^2\theta-4\cos^2\theta = -2\cos^2\theta\leq 0.
\end{align*}

The second case, namely that of $\langle w,\nu_x\rangle\leq 0$, is proved similarly; this time we define $\alpha\in[0,\pi/2)$ by $\langle w,\nu_x\rangle=-\sin\alpha$, directing $\partial_n^x$ such that $\langle w,\partial^x_n\rangle=\cos\alpha$. In this case we define $\theta\in[0,\pi/2)$ and the orientation of $\partial^y_n$ to satisfy the conditions $\langle \partial^y_n,\partial^x_n\rangle=-\cos 2\theta$ and $\langle \partial^y_n,\nu_x\rangle=\sin 2\theta$. A similar calculation as in the first case then yields $\theta=\alpha$ but \eqref{eq:ynidentity} instead implies $\langle \partial^y_n,w\rangle=-\cos\theta$. We now compute:
\begin{align*}
\delta_{qj}-\langle\partial^y_j,\partial^x_q\rangle +2\langle w,\partial^x_q\rangle\langle w,\partial^y_j-\partial^x_j\rangle& =1+\cos(2\theta)-2\cos\theta(\cos\theta+\cos\theta)\\
&= 2\cos^2\theta-4\cos^2\theta = -2\cos^2\theta\leq 0.
\end{align*}

\end{proof}

The matrix $\dot F_x^{ij}\left(Z\delta_{ip}-\xh_{ip}\right)\delta^{pq}$ is non-negative definite and symmetric (since the factors
are each positive definite and commute), so in particular the component with $j=q=n$ is non-negative.  We therefore conclude that
$$
0\leq -\frac{\partial\phi}{\partial t}+\dot F_x^{ij}\nabla_i\nabla_j\phi+\phi\dot F_x^{ij}\xh_{ip}g_x^{pq}\xh_{qj},
$$
which completes the proof that $\overline Z$ is a viscosity subsolution of \eqref{eq:vareqn}.
In the case where $F$ is convex and we consider $\underline{Z}$ instead of $\overline Z$, all inequalities are reversed and we deduce that $\underline{Z}$ is a supersolution of \eqref{eq:vareqn}.
\end{proof}

\section{Conclusions and remarks}

We mention here some immediate implications of the non-collapsing result:
\begin{enumerate}[label={(\arabic*).}]
\item Interior non-collapsing for concave $F$ rules out blow-up limits such as the product of the grim reaper with $\RR^{n-1}$ (if the initial hypersurface has positive $F$), since this has the interior sphere curvature $\overline Z$ asymptotically constant while the speed $F$ approaches zero, violating Corollary \ref{cor:collapse-scale}.  The exterior non-collapsing does not appear to rule out this possibility.  Note that without the assumption of embeddedness, such singularities do indeed occur, even in mean curvature flow.
\item In the case of mean curvature flow where both interior and exterior non-collapsing hold, we are able to deduce directly that for uniformly convex hypersurfaces all principal curvatures are comparable, implying a simple proof of the Huisken and Gage-Hamilton theorems on the asymptotic behaviour for convex solutions \cites{Huiskenconvex, GH}.  If only one-sided non-collapsing holds then we cannot immediately conclude such a strong result, but nevertheless the convergence arguments in the convex case become rather easy:   For example, in the case where $F$ is convex, we have $\underline{Z}(x,t)\geq \varepsilon F(x,t)\geq \varepsilon \kappa_{\max}(x,t)$, from which it follows that the circumradius (bounded by the reciprocal of $\underline{Z}(x,t)$ for any $x$) is bounded by $\varepsilon^{-1}$ times the inradius.   No such result holds in the case where $F$ is concave, however --- this should not be surprising since there are examples of concave, homogeneous degree one functions $F$ such that convex hypersurfaces can evolve to be non-convex under equation \eqref{E:theflow} (see \cite{AMZConvexHypersurfaces}*{Example 1}).
\item As in the case of mean curvature flow, analogues of Corollary \ref{cor:collapse-scale} hold with $F$ replaced by any positive solution of the linearized flow \eqref{eq:vareqn}.  In particular we can allow star-shaped initial hypersurfaces even if $F$ is not positive, by using the solution
$\langle X,\nu\rangle +2tF$ of \eqref{eq:vareqn}.
\end{enumerate}

\end{document}